\documentclass[english]{amsart}
\usepackage{amsmath}
\usepackage{amssymb}
\usepackage{amsthm}
\usepackage{babel}
\usepackage{graphicx}

\newtheorem*{prop}{Proposition}
\newtheorem*{teor}{Theorem}
\newtheorem*{coro}{Corollary}
\newtheorem*{lem}{Lemma}
\theoremstyle{definition}

\newtheorem*{exe}{Exercise}
\newtheorem*{rem}{Remark}
\newtheorem*{defi}{Definition}
\author{Antonio M. Oller Marc\'{e}n}
\title{Origami constructions}
\address{Departamento de Matem\'{a}ticas, Universidad de Zaragoza\\ C/Pedro Cerbuna 12, 50009 Zaragoza (Spain)}
\email{oller@unizar.es}

\date{\today}
\begin{document}
\maketitle 

One of the cornerstones of elementary field theory is the characterization of the real numbers which are constructible using only a ruler (without marks) and a compass. This characterization leads to the impossibility of doubling the cube, trisecting an angle or squaring the circle using only those tools (see chapter 3, section 3 in \cite{ELD} or chapter V, section 1 in \cite{HUN}). Nevertheless, even the ancient Greeks were aware of the possibility of solving those problems using other techniques like mechanical procedures or transcendental curves (see \cite{HEA} for instance). Here we will adopt a different - and somewhat more rigorous - point of view, in particular we will be interested in characterizing the real numbers constructible by paper folding. For slightly different approaches we recommend the reading of \cite{ALP} and chapter 10, section 3 in \cite{COX}. For further relations between origami and mathematics book \cite{HUL} might be worth reading.

Like when working with ruler and compass, we will fix a couple of points in the euclidean plane $\mathbb{R}^2$ which will determine both our origin and `unit segment'. Then the steps that we are allowed to do during any construction are the following:
\begin{itemize}
\item[i)] To draw a line connecting two points already constructed.
\item[ii)] To find the intersection point of two constructible lines.
\item[iii)] To construct the perpendicular bisector of the segment connecting two constructible points.
\item[iv)] To draw the line bisecting any constructed angle.
\item[v)] Given a constructed line $l$ and constructed points $P$ and $Q$, then it is possible to construct a line passing through $Q$ and reflecting $P$ onto $l$.
\item[vi)] Given constructed line $l$ and $m$ and constructed points $P$ and $Q$, then whenever possible, a line reflecting simultaneously $P$ onto $l$ and $Q$ onto $m$ can be constructed.
\end{itemize}
\begin{exe}
\quad 
\begin{enumerate}
\item Show that steps i) to v) can be performed using ruler and compass.
\item Show that the steps allowed in ruler and compass constructions can be performed by paper folding.
\end{enumerate}
\end{exe}

\begin{defi}
An element $x\in\mathbb{R}_{\geq0}$ is said to be \textit{origami constructible} if, starting with the unit segment and using only operations i) to vi), it is possible to construct a segment of length $x$. By definition an element $x\in\mathbb{R}_{<0}$ is \textit{origami constructible} if and only if so is $-x$.
\end{defi}

In what follows we will consider the set $\mathcal{O}=\{ x\in\mathbb{R}\ |\ \textrm{$x$ is origami constructible}\}$. From the previous exercise and the properties of ruler and compass constructions we have the following.

\begin{prop}
$\mathcal{O}$ is a subfield of $\mathbb{R}$ and an extension of $\mathbb{Q}$. Moreover, if $0<a\in\mathcal{O}$, then $\sqrt{a}\in\mathcal{O}$.
\end{prop}

Now we will turn again to step v). Let constructible line $l$ and constructible points $P$ and $Q$ be given, then we consider a parabola with directrix $l$ and focus $P$. We can now construct (see figure 1), if it exists, a tangent to the parabola passing through $Q$. By well-known geometric properties of the parabola, this line reflects $P$ onto $l$ as desired.

\begin{figure}[h]
\center{\includegraphics[width=6 cm]{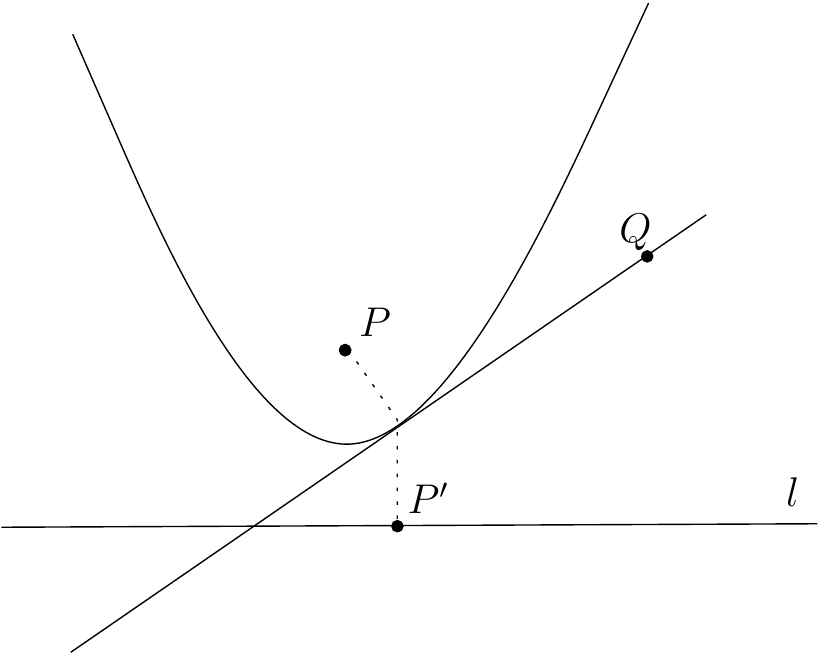}}
\caption{A closer look at step v).}
\end{figure}

Once we have a better understanding of step v) we will look at step vi). It can be seen as two simultaneous applications of step v) thus, if constructible lines $l$ and $m$ and constructible points $P$ and $Q$ are given, step vi) allows us to construct (if it exists, of course) the common tangent to the parabolas with given data as directrices and foci. Let us make a more analytic approach. We will suppose without loss of generality that our parabolas have equations: 
$$\left(y-\frac12 a\right)^2=2bx$$ 
$$y=\frac12 x^2$$
Now consider a simultaneous tangent with slope $\mu$ meeting these parabolas at the respective points $(x_0,y_0)$ and $(x_1,y_1)$. By definition $\displaystyle{\mu=\frac{y_1-y_0}{x_1-x_0}}$, while some easy differentiation yields to
$\displaystyle{\mu=x_1=\frac{b}{y_0-\frac12 a}}$, $\displaystyle{y_1=\frac12 x_1^2}$ and $\displaystyle{x_0=\frac{\left(y_0-\frac12 a\right)^2}{2b}}$. Finally substituting this data in the former expression for $\mu$ one obtains that $$\mu^3+a\mu+b=0.$$
So we have seen that if $a$ and $b$ are origami constructible numbers, so are the real roots of $\mu^3+a\mu+b$. In particular, as the constructibility of a line implies the constructibility of its slope, the following holds.

\begin{prop}
If $a\in\mathcal{O}$ and $\sqrt[3]{a}\in\mathbb{R}$, then $\sqrt[3]{a}\in\mathcal{O}$.
\end{prop}

Now, given a field extension $F/\mathbb{Q}$ with $F\subseteq\mathbb{R}$ we will say that:
\begin{itemize}
\item[i)] A line $l$ belongs to $F$ if it passes through two points with coordinates in $F$. Equivalently if it passes through a point with coordinates in $F$ and its slope is either $\infty$ or it belongs to $F$.
\item[ii)] A parabola $\mathcal{P}$ belongs to $F$ if its focus has coordinates in $F$ and its directrix belongs to $F$.
\end{itemize}

Keeping this definition in mind we can present the following lemma, whose proof will we left as an exercise.

\begin{lem}
Let $F$ be a field with $\mathbb{Q}\subseteq F\subseteq\mathbb{R}$, and let $l$ be a line in $\mathbb{R}^2$ and $\mathcal{P}$ a parabola in $\mathbb{R}^2$. Then:
\begin{enumerate}
\item $l$ belong to $F$ if and only if there are $a,b,c\in F$ such that $l=\{(x,y)\in\mathbb{R}^2\ |\ ax+by+c=0\}$
\item If $\mathcal{P}$ belongs to $F$, then there are $a,b,c,d,e,f\in F$ with $b^2=4ac$ and $\displaystyle{\det\begin{pmatrix}a & \frac{b}{2} & \frac{d}{2}\\ \frac{b}{2} & c & \frac{e}{2}\\ \frac{d}{2} & \frac{e}{2} & f\end{pmatrix}\neq 0}$ such that $\mathcal{P}=\{(x,y)\in\mathbb{R}^2\ |\ ax^2+bxy+cy^2+dx+ey+f=0\}$. The converse is true if $F$ is closed for square roots.
\end{enumerate}
\end{lem}

Now, we are ready to proof the following proposition.

\begin{prop}
Let $F$ be a field with $\mathbb{Q}\subseteq F\subseteq{R}$. Then:
\begin{enumerate}
\item If $l_1$ and $l_2$ are nonparallel lines that belong to $F$, then the coordinates of $l_1\cap l_2$ are in $F$.
\item If $l_1$ and $l_2$ are lines that belong to $F$, then there exist $u,v\in\mathbb{R}$ with $u^2,v^2\in F$ such that the line bisecting the angle between $l_1$ and $l_2$ belongs to $F(u,v)$.
\item If $P\in\mathbb{R}^2$ is a point with coordinates in $F$ and $\mathcal{P}$ is a parabola that belongs to $F$, then there exists $u\in\mathbb{R}$ with $u^2\in F$ such that the lines tangent to $\mathcal{P}$ passing through $P$ belong to $F(u)$.
\item If $\mathcal{P}_1$ and $\mathcal{P}_2$ are parabolas that belong to $F$, then there exist $u,v\in\mathbb{C}$ with $u^2\in F$ and $v^3\in F(u)$ such that the common tangent (if it exists) to $\mathcal{P}_1$ and $\mathcal{P}_2$ belongs to $F(u,v)\cap\mathbb{R}$.
\end{enumerate}
\end{prop}
\begin{proof}
Just some boring computations together with the previous lemma and the form of the solutions to quadratic and cubic equations. Details will be left to the reader.
\end{proof}

\begin{teor}
Let $a\in\mathbb{R}$, then $a\in\mathcal{O}$ if and only if there are $0\leq n\in\mathbb{Z}$ and subfields $F_0,F_1,\dots F_n$ of $\mathbb{R}$ such that $\mathbb{Q}=F_0\leq F_1\leq\dots\leq F_n$, $[F_i:F_{i-1}]=2\ \textrm{or}\ 3$ for any $1\leq i\leq n$ and $a\in F_n$. In particular, if $a\in\mathcal{O}$, then $[\mathbb{Q}(a):\mathbb{Q}]=2^r3^s$ for some $0\leq r,s\in\mathbb{Z}$.
\end{teor}
\begin{proof}
If $a\in\mathcal{O}$, then we can construct it starting from the `unit segment' by successive applications of steps i) to vi). The previous proposition shows that at each step the coordinates of the new point are either in the same field where the coordinates of the previous points are, or in a quadratic extension, or in a quadratic extension of a quadratic extension or in a cubic extension of a quadratic one. In any case the result holds.

Conversely, let us suppose that $F_i=F_{i-1}(u_i)$ with $u_i^2\in F_{i-1}$ and $u_i\notin F_{i-1}$ (quadratic extension) or $u_i^3\in F_{i-1}$ and $u_i,u_i^2\notin F_{i-1}$ (cubic extension). Obviously $F_0=\mathbb{Q}\subseteq\mathcal{O}$, moreover if $F_{i-1}\subseteq\mathcal{O}$ then $u_i^2\in\mathcal{O}$ or $u_i^3\in\mathcal{O}$ and, in any case, as $\mathcal{O}$ is closed for square and cubic roots be have that $u_i\in\mathcal{O}$ and consequently that $F_i\subseteq\mathcal{O}$. Therefore $F_n\subseteq\mathcal{O}$ and since $a\in F_n$ the result holds.

Finally, from $[F_n:\mathbb{Q}(a)][\mathbb{Q}(a):\mathbb{Q}]=[F_n:\mathbb{Q}]=2^r3^s$, the last assertion follows.
\end{proof}

\begin{rem}
As in the case of ruler and compass constructions we could prove a converse with some use of `Galois Theory'. It can be proved that a real number $a$ is origami constructible if and only if $[K:\mathbb{Q}]=2^r3^s$ with $K$ being the splitting field of its minimum polynomial over $\mathbb{Q}$.
\end{rem}

\begin{coro}
If $C=\{x\in\mathbb{R}\ |\ x\ \textrm{ruler and compass constructible}\}$, then we have that $\mathbb{Q}\subset C\subset\mathcal{O}\subset\mathbb{R}$.
\end{coro}

\begin{coro}
The classical problems of doubling the cube and trisecting the angle are solvable by paper folding.
\end{coro}
\begin{proof}
It is trivial as they only involve the resolution of cubic equations. 
\end{proof}

The preceding corollary shows the possibility of solving those problems. It states the existence of a paper folding construction which calculates $\sqrt[3]{2}$ or that trisects a given angle but, as it often happens in mathematics, it does not give any clue about how to find such a construction. Just to convince ourselves we will present now those constructions.

\begin{figure}[h]
\center{\includegraphics[width=8.5 cm]{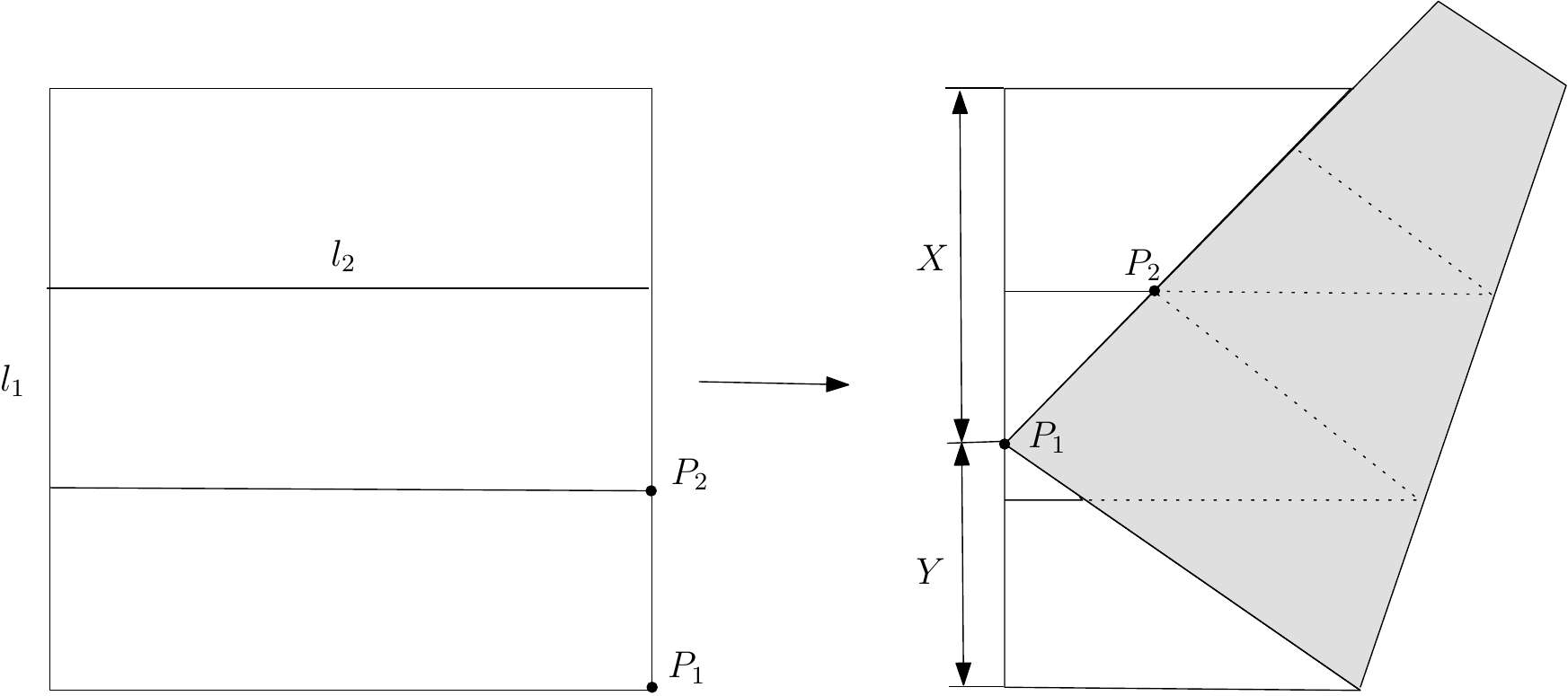}}
\caption{Doubling the cube.}
\end{figure}

\begin{figure}[h]
\center{\includegraphics[width=8.5 cm]{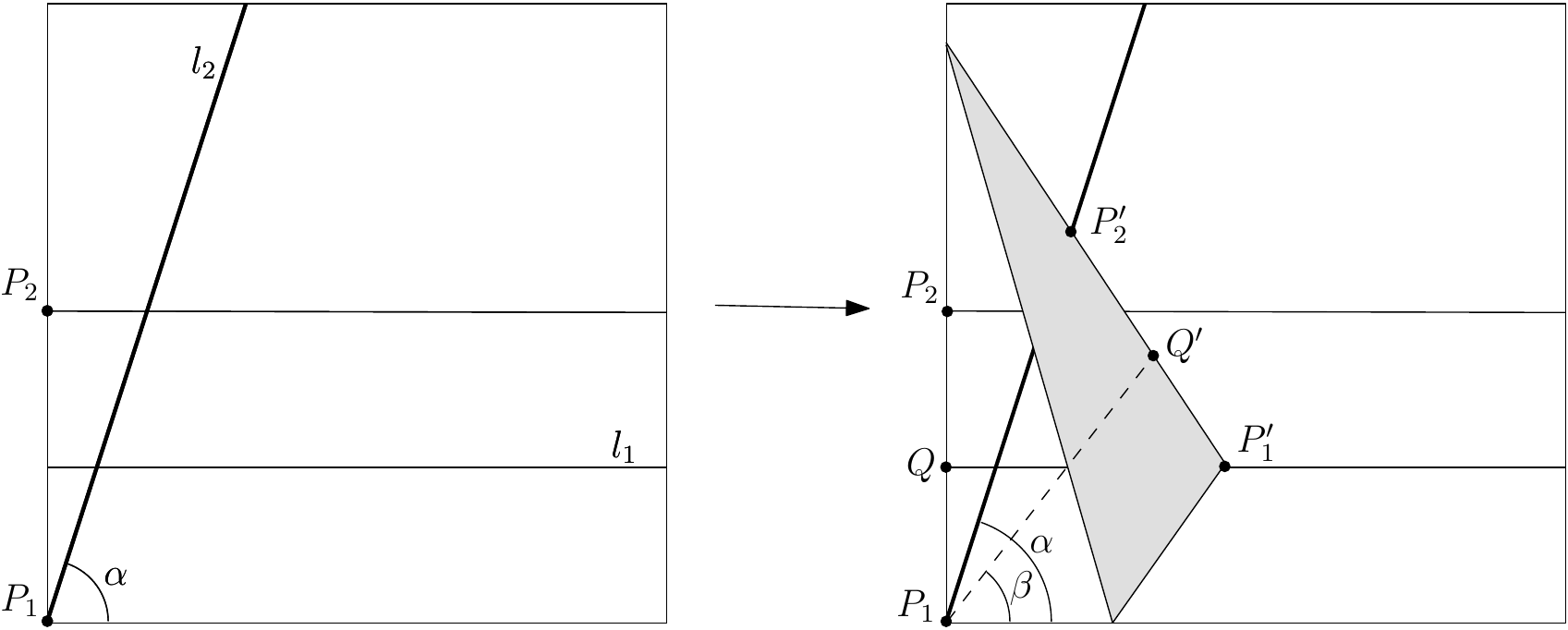}}
\caption{Trisecting an angle.}
\end{figure}

Take a square of paper and divide it horizontally into three equal parts. Fold it in such a way that point $P_1$ reflects onto line $l_1$ and $P_2$ onto $l_2$ (see figure 2).

\begin{exe}
Show that $\displaystyle{\frac{X}{Y}=\sqrt[3]{2}}$.
\end{exe}

Now, take a square of paper and put any angle $\alpha$ between $\displaystyle{\frac{\pi}{4}}$ and $\displaystyle{\frac{\pi}{2}}$ in the bottom left corner of the paper. Divide the paper horizontally into two parts and the divide the half below again into two parts. Then, fold the paper in such a way that $P_1$ reflects onto $l_1$ and $P_2$ onto $l_2$ (see figure 3).

\begin{exe}
Show\footnote{As somebody said once, ``Geometry is the art of correct reasoning on incorrect figures'', so be careful and do not trust the picture.} that $\beta=\displaystyle{\frac23\alpha}$.
\end{exe}

Finally we will characterize the regular polygons which are constructible by paper folding. We have the following. 

\begin{prop}
The regular polygon of $n$ sides is constructible by paper folding if and only if $n=2^r3^sp_1\dots p_t$ where $0\leq r,s,t\in\mathbb{Z}$ and $p_1,\dots, p_t$ are distinct primes of the form $p_i=2^{a_i}3^{b_i}+1$.
\end{prop}
\begin{proof}
Let us denote the regular polygon of $n$ sides by $P(n)$. Then, it is clear that $P(n)$ is constructible by paper folding if and only if $\displaystyle{\cos\frac{2\pi}{n}\in\mathcal{O}}$.

First note that if $m|n$ and $P(n)$ is constructible, then so is $P(m)$. On the other hand, if $(m,n)=1$ and both $P(n)$ and $P(m)$ are constructible, then so is $P(mn)$. In fact, it can be seen (we leave some details as an exercise) that $\mathbb{Q}\left(\displaystyle{\cos\frac{2\pi}{mn}}\right)=\mathbb{Q}\left(\displaystyle{\cos\frac{2\pi}{m},\cos\frac{2\pi}{n}}\right)\subseteq\mathcal{O}$. These considerations allow us to reduce to the primary case.

Obviously $P(2^r)$ and $P(3^s)$ are constructible for every $0\leq r,s\in\mathbb{Z}$, so let $p\neq 2,3$ be a prime. If we put $z=e^{\frac{2\pi i}{p}}$, then $\displaystyle{\cos\frac{2\pi}{p}=\frac{z+\overline{z}}{2}}$. As $[\mathbb{Q}(z):\mathbb{Q}(z+\overline{z})]=2$ and $[\mathbb{Q}(z):\mathbb{Q}]=p-1$ then $[\mathbb{Q}(z+\overline{z}):\mathbb{Q}]=\displaystyle{\frac{p-1}{2}}$. Thus, $\displaystyle{\cos\frac{2\pi}{p}\in\mathcal{O}}$ if and only if $\displaystyle{\frac{p-1}{2}}=2^{a}3^b$ and the proof is complete.
\end{proof}

\end{document}